# Solutions around a Regular $\alpha$ Singular Point of a Sequential Conformable Fractional Differential Equation

Emrah Ünal[*], Ahmet Gökdoğan[**], Ercan Çelik[***],

[*] Department of Elementary Mathematics Education, Artvin Çoruh University, 08100 Artvin, Turkey
emrah.unal@artvin.edu.tr
[**] Department of Mathematical Engineering, Gümüşhane University, 29100 Gümüşhane, Turkey,
gokdogan@gumushane.edu.tr
[***] Department of Mathematics, Atatürk University, 25100 Erzurum, Turkey,
ecelik@atauni.edu.tr

**Abstract**

In this work, the power series solutions are given around a regular-singular point, in the case of variable coefficients for homogeneous sequential linear conformable fractional differential equations of order $2\alpha$.

**Key words** conformable fractional derivative, regular $\alpha$ singular point, sequential conformable fractional differential equation,

## 1. Introduction

Though fractional derivative idea is more than 300 years old, intensive studies about fractional calculus were carried out many researches in the last and present centuries. Several mathematicians as Liouville, Riemann, Weyl, Fourier, Abel, Leibniz, Grunwald and Letnikov made major contributions to the theory of fractional calculus. The most popular ones of fractional derivative definitions are Grunwald-Letnikov, Riemann-Liouville and Caputo definitions. For Grunwald-Letnikov, Riemann-Liouville, Caputo and other definitions and the characteristics of these definitions, we refer to reader to [1-3].

(I) Grunwald-Letnikov definition:
$$_aD_x^\alpha f(x) = \lim_{h \to 0} h^{-\alpha} \sum_{j \to 0}^{\frac{x-a}{h}} (-1)^j \binom{\alpha}{j} f(x - jh)$$

(II) Riemann Liouville definition:
$$D_x^\alpha f(x) = \frac{1}{\Gamma(n-\alpha)} \left(\frac{d}{dx}\right)^n \int_0^x (x-t)^{n-\alpha-1} f(t) dt, \quad n - 1 < \alpha \leq n$$

(III) Caputo definition:
$$D_x^\alpha f(x) = \frac{1}{\Gamma(n-\alpha)} \int_0^x (x-t)^{n-\alpha-1} \left(\frac{d}{dx}\right)^n f(t) dt, \quad n - 1 < \alpha \leq n$$

Recently, Khalil at al. give a new definition of fractional derivative and fractional integral [4]. This new definition benefit from a limit form as in usual derivatives. They also proved the product rule, the fractional Rolle's theorem and mean value theorem. This new theory is improved by Abdeljawad. For instance, he gives Taylor power series representation and Laplace transform of certain functions, fractional integration by parts formulas, chain rule and Gronwall inequality [5].

The existence of solutions around an ordinary point of conformable fractional differential equation of order $2\alpha$, the solution of conformable fractional Hermite differential equation and fractional Hermite polynomials is analyzed in [6]. Legendre conformable fractional equation and Legendre fractional polynomials are studied in [7].

In this work, we analyze the existence of solutions around a regular $\alpha$ singular point of conformable fractional differential equation of order $2\alpha$.

## 2. Conformable Fractional Calculus

**Definition 2.1.** [5] Given a function $f:[a,\infty) \to \mathbb{R}$. Then the left conformable fractional derivative of $f$ order $\alpha$ is defined by

$$(T_\alpha^a f)(x) = \lim_{\varepsilon \to 0} \frac{f(x + \varepsilon(x-a)^{1-\alpha}) - f(x)}{\varepsilon}$$

for all $x > a, \alpha \in (0,1]$. When $a = 0$, it is written as $T_\alpha$. If $(T_\alpha f)(x)$ exists on $(a,b)$ then $(T_\alpha^a f)(a) = \lim_{x \to a^+} (T_\alpha^a f)(x)$.

**Definition 2.2.** [5] Given a function $f:(-\infty, b] \to \mathbb{R}$. Then the right conformable fractional derivative of $f$ order $\alpha$ is defined by

$$({}_\alpha^b T f)(x) = -\lim_{\varepsilon \to 0} \frac{f(x + \varepsilon(b-x)^{1-\alpha}) - f(x)}{\varepsilon}$$

for all $x < b, \alpha \in (0,1]$. If $({}_\alpha T f)(x)$ exists on $(a,b)$ then $({}_\alpha^b T f)(b) = \lim_{x \to b^-} ({}_\alpha^b T f)(x)$.

**Theorem 2.1.** [4] Let $\alpha \in (0,1]$ and $f, g$ be $\alpha$-differentiable at a point $x > 0$. Then

(1) $\frac{d^\alpha}{dx^\alpha}(af + bg) = a\frac{d^\alpha f}{dx^\alpha} + b\frac{d^\alpha g}{dx^\alpha}$, for all $a, b \in \mathbb{R}$

(2) $\frac{d^\alpha}{dx^\alpha}(x^p) = px^{p-\alpha}$, for all $p \in \mathbb{R}$

(3) $\frac{d^\alpha}{dx^\alpha}(\lambda) = 0$, for all constant functions $f(x) = \lambda$

(4) $\frac{d^\alpha}{dx^\alpha}(fg) = f\frac{d^\alpha}{dx^\alpha}(g) + g\frac{d^\alpha}{dx^\alpha}(f)$

(5) $\frac{d^\alpha}{dx^\alpha}(f/g) = \frac{g\frac{d^\alpha}{dx^\alpha}(f) - f\frac{d^\alpha}{dx^\alpha}(g)}{g^2}$

(6) If, in addition, $f$ is differentiable, then $\frac{d^\alpha}{dx^\alpha}(f(x)) = x^{1-\alpha}\frac{df}{dx}(x)$.

**Theorem 2.2.** [5] Assume $f$ is infinitely $\alpha$-differentiable function, for some $0 < \alpha \leq 1$ at a neighborhood of a point $x_0$. Then $f$ has the fractional power series expansion:

$$f(x) = \sum_{k=0}^{\infty} \frac{\left(^{(k)}T_\alpha^{x_0}f\right)(x_0)(x-x_0)^{k\alpha}}{\alpha^k k!}, \qquad x_0 < x < x_0 + R^{1/\alpha}, \qquad R > 0.$$

Here, $\left(^{(k)}T_\alpha^{x_0}f\right)(x_0)$ means the application of the fractional derivative $k$ times.

## 3. Conformable Fractional Differential Equation and Solutions around a regular $\alpha$ singular point

The most general sequential linear homogeneous (left) conformable fractional differential equation is

$$^{(n)}T_\alpha^a y + a_{n-1}(x)\,^{(n-1)}T_\alpha^a y + \ldots + a_1(x)T_\alpha^a y + a_0(x)y = 0, \qquad (1)$$

where $^{(n)}T_\alpha^a y = T_\alpha^a T_\alpha^a \ldots T_\alpha^a y$, n times.

**Definition 3.1.** Let $\alpha \in (0,1]$, $x_0 \in [a,b]$, $N(x_0)$ be a neighborhood of $x_0$ and $f(x)$ be a real function defined on $[a,b]$. In this case $f(x)$ is said to be $\alpha$-analytic at $x_0$ if $f(x)$ can be expressed as a series of natural powers of $(x - x_0)^\alpha$ for all $x \in N(x_0)$. In other word, $f(x)$ can be expressed as following:

$$\sum_{k=0}^{\infty} c_k (x - x_0)^{k\alpha} \qquad (c_k \in R)$$

This series being definitely convergent for $|x - x_0| < \delta$ ($\delta > 0$). $\delta$ is the radius of convergence of the series.

**Definition 3.2.** Let $\alpha \in (0,1]$, $x_0 \in [a,b]$ and the functions $a_k(x)$ be $\alpha$-analytic at $x_0 \in [a,b]$ for $k = 0,1,2,\ldots,n-1$. In this case, the point $x_0 \in [a,b]$ is said to be an $\alpha$-ordinary point of equation (1). If a point $x_0 \in [a,b]$ is not $\alpha$-ordinary point, then it is said to be $\alpha$ singular.

**Definition 3.3.** Let $\alpha \in (0,1]$, $a_k(x)$ be $\alpha$-singular at $x_0 \in [a,b]$ for $k = 0,1,2,\ldots,n-1$. If the functions $(x - x_0)^{(n-k)\alpha} a_k(x)$ are $\alpha$-analytic at the point $x_0 \in [a,b]$ for $k = 0,1,2,\ldots,n-1$, then the point $x_0$ is said to be a regular $\alpha$ singular point of (1). In the contrary case $x_0$ is said to be an essential $\alpha$ singular point.

**Example 3.1.** a) we shall consider following the conformable fractional differential equations

$$x^\alpha T_\alpha y - y = 0$$

$$x^{2\alpha}\,^2T_\alpha y - 2x^\alpha y = 0$$

$$x^{2\alpha}\,^2T_\alpha y - 2x^\alpha T_\alpha y + x^{2\alpha} y = 0$$

the point $x = 0$ is a regular $\alpha$ singular point for the above equations.

b)

$$(x-1)^\alpha T_\alpha y - y = 0$$
$$(x-1)^{2\alpha}\,{}^2T_\alpha y - 2(x-1)^\alpha T_\alpha y + (x-1)^{2\alpha} y = 0$$

For the these equation, the point $x = 1$ is a regular $\alpha$ singular point.

Now, we consider the following homogeneous sequential linear fractional differential equation of order $2\alpha$:

$$(x-x_0)^{2\alpha} T_\alpha^{x_0} T_\alpha^{x_0} y + (x-x_0)^\alpha p(x) T_\alpha^{x_0} y + q(x) y = 0 \qquad (2)$$

where $\alpha \in (0,1]$. If the point $x_0$ is a regular $\alpha$ singular point of equation (2), then this point is $\alpha$-analytic point for functions $p(x)$ and $q(x)$. In this case, functions $p(x)$ and $q(x)$, respectively, have the following series expansions:

$$p(x) = \sum_{k=0}^{\infty} p_k (x-x_0)^{k\alpha} \qquad (0 < x - x_0 < \delta_1;\; \delta_1 > 0)$$

and

$$q(x) = \sum_{k=0}^{\infty} q_k (x-x_0)^{k\alpha} \qquad (0 < x - x_0 < \delta_1;\; \delta_1 > 0).$$

For equation (2), suppose we have a solution the form

$$y(x; s) = \sum_{k=0}^{\infty} c_k(s)(x-x_0)^{(k+s)\alpha} \qquad (3)$$

where let be $c_0 \neq 0$, $s$ being a number to be determined.

If we substitute equation (3) and conformable derivatives of equation (3) in the equation (2), then we get

$$c_0 I_0(s)(x-x_0)^{s\alpha} + \sum_{k=1}^{\infty}\left[c_k I_0(k+s) + \sum_{j=0}^{k-1} c_j I_{k-j}(j+s)\right](x-x_0)^{(k+s)\alpha} = 0$$

where

$$I_0(s) = \alpha^2 s(s-1) + \alpha s p_0 + q_0 \qquad (4)$$
$$I_m(s) = p_m \alpha s + q_m \qquad (5)$$

Equation (4) is called fractional indicial equation of equation (2). The coefficients $c_k$ is

$$c_k = -\frac{\sum_{j=0}^{k-1} a_j I_{k-j}(j+s)}{I_0(k+s)} \qquad (6)$$

**Theorem 3.1.** Let $\alpha \in (0,1]$ and $x_0$ be a regular $\alpha$ singular point of the equation

$$(x-x_0)^{2\alpha} T_\alpha^{x_0} T_\alpha^{x_0} y + (x-x_0)^\alpha p(x) T_\alpha^{x_0} y + q(x) y = 0.$$

Let $s_1, s_2$ be discint and $s_1 - s_2 \neq n$ for $n \in N$ two real roots of the fractional indicial equation. Then, there exists two linearly independent solution to the equation (2) as following:

$$y_1(x; s_1) = \sum_{k=0}^{\infty} c_k(s_1)(x - x_0)^{(k+s_1)\alpha} \tag{7}$$

$$y_2(x; s_2) = \sum_{k=0}^{\infty} c_k(s_2)(x - x_0)^{(k+s_2)\alpha} \tag{8}$$

for $x \in (x_0, x_0 + \rho)$ with $\rho = \min\{\delta_1, \delta_2\}$ and initial conditions $c_0 = y(x_0)$, $\alpha c_1 = T_\alpha y(x_0)$. Since $x_0$ is a regular $\alpha$ singular point of equation (2), by Definition 3.1 and Definition 3.3, it can be denoted that,

$$p(x) = \sum_{k=0}^{\infty} p_k(x - x_0)^{k\alpha} \quad (x \in [x_0, x_0 + \delta_1]; \delta_1 > 0) \tag{9}$$

and

$$q(x) = \sum_{k=0}^{\infty} q_k(x - x_0)^{k\alpha} \quad (x \in [x_0, x_0 + \delta_2]; \delta_2 > 0). \tag{10}$$

**Proof.** We must prove that series equation (3) converges for $\in (x_0, x_0 + \rho)$. Let be $s = s_1$ and $s = s_2$ such that $s_1 - s_2$ is not a positive integer. We note that

$$I_0(s) = \alpha^2(s - s_1)(s - s_2).$$

Hence, the following equations can be written,

$$I_0(s_1 + k) = \alpha^2 k(k + s_1 - s_2),$$

$$I_0(s_2 + k) = \alpha^2 k(k + s_2 - s_1).$$

Therefore, we get

$$I_0(s_1 + k) \geq \alpha^2 k(k - |s_1 - s_2|) \tag{11}$$

$$I_0(s_2 + k) \geq \alpha^2 k(k - |s_2 - s_1|) \tag{12}$$

Now, let be any number such that $0 < r < \rho$. Series in (9) and (10) converge for $x \in [x_0, x_0 + r]$. Hence, there is a constant number $M > 0$ such that

$$|p_j| r^{j\alpha} \leq M \quad (j \in N) \tag{13}$$

$$|q_j| r^{j\alpha} \leq M \quad (j \in N) \tag{14}$$

Using (11), (12), (13) and (14) in (6), we have

$$|c_k(s_1)| \leq \frac{M}{r^{k\alpha}} \frac{\sum_{j=0}^{k-1} \alpha(j+1+|s_1|) r^{j\alpha} |c_j(s_1)|}{\alpha^2 k(k - |s_1 - s_2|)} \tag{15}$$

Now, let be an integer number $N$ such that

$$N - 1 \leq |s_1 - s_2| < N.$$

We define

$$C_0 = c_0(s_1) = 1, C_1 = |c_1(s_1)|, \ldots, C_{N-1} = |c_{N-1}(s_1)|$$

and let the coefficients $C_k$ for $k \geq N$ be defined by

$$C_k = \frac{M}{r^{k\alpha}} \frac{\sum_{j=0}^{k-1} \alpha(j+1+|s_1|) r^{j\alpha} C_j(s_1)}{\alpha^2 k(k - |s_1 - s_2|)}. \tag{16}$$

From the definition of $C_k$ and equation (15), we see that

$$|c_k(s_1)| \leq C_k \quad k = 0,1,2,\ldots$$

We prove that the series

$$\sum_{k=0}^{\infty} C_k(x-x_0)^{k\alpha} \tag{17}$$

is convergent for $x \in (x_0, x_0 + \rho)$. By using (16), we obtain that

$$r^\alpha \alpha^2(k+1)(k+1-|s_1-s_2|)C_{k+1} = \alpha^2(k)(k-|s_1-s_2|)C_k + \alpha M(k+1+|s_1|)C_k.$$

Hence,

$$\frac{C_{k+1}}{C_k} = \frac{\alpha^2(k)(k-|s_1-s_2|) + \alpha M(k+1+|s_1|)}{r^\alpha \alpha^2(k+1)(k+1-|s_1-s_2|)}$$

is obtained. By the help of the ratio test, we have that

$$\lim_{k \to \infty} \left| \frac{C_{k+1}(x-x_0)^{(k+1)\alpha}}{C_k(x-x_0)^{k\alpha}} \right| = \left(\frac{|x-x_0|}{r}\right)^\alpha < 1.$$

Thus, the series (17) converges for $x \in [x_0, x_0 + r]$. This implies that the series (3) converges for $x \in [x_0, x_0 + r]$ and $s_1$. Since $r$ was any number satisfying $0 < r < \rho$, the series (3) converges for $x \in (x_0, x_0 + \rho)$.

Similarly, the same computations with $s_1$ replaced by $s_2$ everywhere show that the series (3) converges for $x \in (x_0, x_0 + \rho)$ and $s_2$.

**Theorem 3.2.** Let $x_0$ be a regular $\alpha$-singular point of equation (2). For this equation, $p(x)$ and $q(x)$ have fractional power series expansion, for $x \in (x_0, x_0 + \rho)$ with $\rho > 0$. Let $s_1$ and $s_2$ be two real roots of the fractional indicial equation (4).

If $s_1 = s_2$, there are two linearly independent solution and these solutions have, respectively, the following forms:

$$y_1(x; s_1) = \sum_{k=0}^{\infty} c_k(x-x_0)^{(k+s_1)\alpha} \quad (c_0 \neq 0)$$

$$y_2(x; s_1) = \ln|x-x_0|y_1(x; s_1) + \sum_{k=0}^{\infty} b_k(x-x_0)^{(k+s_1+1)\alpha}$$

for $x \in (x_0, x_0 + \rho)$.

If $s_1 - s_2$ is a positive integer, linearly independent solutions have, respectively, the following forms:

$$y_1(x; s_1) = \sum_{k=0}^{\infty} c_k(x-x_0)^{(k+s_1)\alpha} \quad (c_0 \neq 0)$$

$$y_2(x; s_1) = C.\ln|x-x_0|y_1(x; s_1) + \sum_{k=0}^{\infty} b_k(x-x_0)^{(k+s_2)\alpha}$$

for $x \in (x_0, x_0 + \rho)$ where $C$ is a constant and it may happen zero.

**Proof:** For $s_1 \geq s_2$, according to Theorem 3.1.

$$y_1(x; s_1) = \sum_{k=0}^{\infty} c_k(s_1)(x-x_0)^{(k+s_1)\alpha} \tag{18}$$

is first solution of (2). We rewrite equation (2) as following

$$T_\alpha^{x_0} T_\alpha^{x_0} y + P(x) T_\alpha^{x_0} y + Q(x) y = 0$$

where $P(x) = \frac{p(x)}{(x-x_0)^\alpha}$, $Q(x) = \frac{q(x)}{(x-x_0)^{2\alpha}}$.

By the help of conformable Abel's formulas in [8], we write

$$y_2(x) = y_1(x) I_\alpha^{x_0}\left(\frac{e^{-I_\alpha^{x_0}(P(x))}}{[y_1(x)]^2}\right). \tag{19}$$

Now, let be $s_1 - s_2 = N$ such that $N$ is non-negative integer. Hence, $s_1$ and $s_2 = s_1 - N$ are roots of indicial equation. Therefore, we get

$$-p_0 - 2\alpha s_1 = \alpha(-1 - N) \tag{20}$$

$$P(x) = \frac{p_0 + p_1(x-x_0)^\alpha + p_2(x-x_0)^{2\alpha} + \cdots}{(x-x_0)^\alpha} = \frac{p_0}{(x-x_0)^\alpha} + p_1 + p_2(x-x_0)^\alpha + \cdots$$

Hence,

$$exp\left(-I_\alpha^{x_0}(P(x))\right) = exp\left(-I_\alpha^{x_0}\left(\frac{p_0}{(x-x_0)^\alpha} + p_1 + p_2(x-x_0)^\alpha + \cdots\right)\right)$$

$$= exp\left(-p_0 ln|x-x_0| - \frac{p_1}{\alpha}(x-x_0)^\alpha - \frac{p_2}{2\alpha}(x-x_0)^{2\alpha} - \cdots\right)$$

$$= (x-x_0)^{-p_0} exp\left(-\frac{p_1}{\alpha}(x-x_0)^\alpha - \frac{p_2}{2\alpha}(x-x_0)^{2\alpha} - \cdots\right).$$

That is,

$$exp\left(-I_\alpha^{x_0}(P(x))\right) = (x-x_0)^{-p_0}(1 + A_1(x-x_0)^\alpha + A_2(x-x_0)^{2\alpha} + \cdots) \tag{21}$$

Now, let us choose $c_0 = 1$ and we substitute (18) and (21) in (19). In this case, we get

$$y_2(x) = y_1(x) I_\alpha^{x_0}\left(\frac{(x-x_0)^{-p_0}(1 + A_1(x-x_0)^\alpha + A_2(x-x_0)^{2\alpha} + \cdots)}{(x-x_0)^{2s_1\alpha}(1 + c_1(x-x_0)^\alpha + c_2(x-x_0)^{2\alpha} + \cdots)^2}\right).$$

$$y_2(x) = y_1(x) I_\alpha^{x_0}\left(\frac{(x-x_0)^{-p_0-2s_1\alpha}(1 + A_1(x-x_0)^\alpha + A_2(x-x_0)^{2\alpha} + \cdots)}{(1 + B_1(x-x_0)^\alpha + B_2(x-x_0)^{2\alpha} + \cdots)}\right)$$

$$= y_1(x) I_\alpha^{x_0}\left((x-x_0)^{\alpha(-1-N)}(1 + C_1(x-x_0)^\alpha + C_2(x-x_0)^{2\alpha} + \cdots)\right).$$

For $N = 0$, that is $s_1 = s_2$, we have

$$y_2(x) = y_1(x) I_\alpha^{x_0}((x-x_0)^{-\alpha} + C_1 + C_2(x-x_0)^\alpha + \cdots)$$

$$= y_1(x) ln(x-x_0) + y_1(x)\left(\frac{C_1}{\alpha}(x-x_0)^\alpha + \frac{C_2}{2\alpha}(x-x_0)^{2\alpha} + \cdots\right)$$

$$= y_1(x) ln(x-x_0) + (x-x_0)^{s_1\alpha}(1 + c_1(x-x_0)^\alpha +$$

$$\cdots)\left(\frac{C_1}{\alpha}(x-x_0)^\alpha + \frac{C_2}{2\alpha}(x-x_0)^{2\alpha} + \cdots\right)$$

$$= y_1(x)\ln(x-x_0) + (x-x_0)^{s_1\alpha}(b_0(x-x_0)^\alpha + b_1(x-x_0)^{2\alpha} + b_3(x-x_0)^{3\alpha} \ldots)$$

Consequently, for $s_1 = s_2$, the general form of second solution is

$$y_2(x) = y_1(x)\ln(x-x_0) + (x-x_0)^{(s_1+1)\alpha} \sum_{k=0}^{\infty} b_k(x-x_0)^{k\alpha}.$$

For $N > 0$, that is $s_1 - s_2 = N$, we have

$$y_2(x) = y_1(x)I_\alpha^{x_0}\left((x-x_0)^{\alpha(-1-N)}(1 + C_1(x-x_0)^\alpha + C_2(x-x_0)^{2\alpha} + \cdots + C_N(x-x_0)^{N\alpha} + \cdots)\right)$$

$$= y_1(x)I_\alpha^{x_0}\left(\left(\frac{C_N}{(x-x_0)^\alpha} + \frac{1}{(x-x_0)^{(1+N)\alpha}} + \frac{C_1}{(x-x_0)^{N\alpha}} + \cdots\right)\right)$$

$$= C_N y_1(x)\ln(x-x_0) + y_1(x)\left(\frac{(x-x_0)^{-N\alpha}}{-N\alpha} + \frac{C_1(x-x_0)^{(-N+1)\alpha}}{(-N+1)\alpha} + \cdots\right)$$

$$= C_N y_1(x)\ln(x-x_0) + (x-x_0)^{(s_2+N)\alpha}(1 + c_1(x-x_0)^\alpha + \cdots)(x-x_0)^{-N\alpha}\left(-\frac{1}{N\alpha} + \frac{c_1(x-x_0)^\alpha}{(-N+1)\alpha} + \cdots\right).$$

Hence, for $N > 0$, the general form of second solution is

$$y_2(x) = C_N y_1(x)\ln(x-x_0) + y_1 + (x-x_0)^{s_2\alpha} \sum_{k=0}^{\infty} b_k(x-x_0)^{k\alpha}.$$

where $b_0 = -\frac{c_0}{N} \neq 0$.

## 4. Conclusion

In this work, power series solutions are given around an regular $\alpha$ singular point in homogeneus sequential linear differential equation of conformable fractional of order $2\alpha$ with variable coefficients. Firstly, definitions of $\alpha$-ordinary point and regular $\alpha$ singular point is given. Then, for distinct roots of indicial equation which the differences between them is not positive integer, general form of solutions is given. Finally, for equal roots and distinct roots which the differences between them is integer, general form of solutions is obtained. It is appeared that the results obtained in this work correspond to results which are obtained in ordinary case.